\documentclass{amsart}

\newtheorem{theorem}{Theorem}[section]
\newtheorem{lemma}[theorem]{Lemma}
\newtheorem{corollary}[theorem]{Corollary}

\theoremstyle{definition}
\newtheorem{definition}[theorem]{Definition}
\newtheorem{example}[theorem]{Example}

\newtheorem{proposition}[theorem]{Proposition}
\newtheorem{remark}[theorem]{Remark}

\numberwithin{equation}{section}

\usepackage{geometry}
\begin{document}

\title[]{Quadratic Chern Forms and Rigidity of \\Hermitian Metrics}

\author[]{Jieming Yang}

\address{Department of mathematics, Wenzhou University, Zhejiang, China.}
\email{20230120@wzu.edu.cn}



\thanks{The author thanks Professor Jixiang Fu. The author also thanks ChatGPT for conversations, and the calculation checks, and the English writing. The author takes full responsibility for the manuscript.}

\date{}

\begin{abstract}
	We relate the Miyaoka-Yau density of a Hermitian metric to the Chern form's bound for the Calabi functional and prove a pointwise upper bound with optimal torsion coefficient in complex dimension at least three. For connected balanced manifolds satisfying the weak second Chern-Einstein equation, the density is nonnegative and vanishes identically precisely for Chern-flat metrics and K\"ahler metrics of nonzero constant holomorphic sectional curvature. On a compact threefold admitting a pluriclosed metric, a single integral of the difference of the first and second Chern-Ricci forms detects whether a given balanced metric is K\"ahler.
\end{abstract}

\maketitle

\section{Introduction}
Let $X$ be a complex manifold with $\dim_{\mathbb{C}}X=n\geq2$, equipped with a Hermitian metric $\omega$. If $\partial\bar{\partial}\,\omega^{n-1}=0$, $\omega$ is called Gauduchon; if $\mathrm{d}\,\omega^{n-1}=0$, $\omega$ is called balanced; if $\partial\bar{\partial}\,\omega^{n-2}=0$, $\omega$ is called astheno-K\"{a}hler; if $\mathrm{d}\,\omega=0$, $\omega$ is called K\"{a}hler.

Let $\nabla$ be the Chern connection of $\omega$, with curvature $R_\nabla$, torsion $T$ and Gauduchon torsion (1,0)-form $\tau$. Write $\rho_\nabla$, $\kappa_\nabla\in\mathcal{A}^{1,1}_{\mathbb{R}}(X)$ for the first and second Chern-Ricci forms, respectively. Their trace-free parts are
\[
\rho^{\circ}_\nabla=\rho_\nabla-\frac{1}{n}s\,\omega, \quad
\kappa^{\circ}_\nabla=\kappa_\nabla-\frac{1}{n}s\,\omega,
\]
where $s=\Lambda_{\omega}\rho_\nabla=\Lambda_{\omega}\kappa_\nabla$ is the first scalar curvature.

Define $\mathcal{D}_\nabla$, $\mathcal{L}_\nabla\in\mathcal{A}^{0}_{\mathbb{R}}(X)$ by
\begin{equation}\label{11233}
	\begin{split}
		\mathcal{D}_\nabla\,\frac{\omega^{n}}{n!}
		&=4\pi^2\left(2(n+1)\,c_2(X, \omega)-n\,c_1(X, \omega)^2\right)\wedge
		\frac{\omega^{n-2}}{(n-2)!},\\
		\mathcal{L}_\nabla\,\frac{\omega^{n}}{n!}
		&=8\pi^2\left(n\,c_1(X, \omega)^2-(n+2)\,c_2(X, \omega)\right)\wedge
		\frac{\omega^{n-2}}{(n-2)!}.
	\end{split}
\end{equation}
The Chern polynomial defining $\mathcal D_\nabla$ is the one appearing in the Miyaoka-Yau inequality. The polynomial defining $\mathcal L_\nabla$ appears in Ping Li's lower bound for the Calabi functional \cite{Li}. For a fixed K\"ahler class
$[\omega]$ on a compact K\"ahler manifold, the bound involves the intersection number
\[
\left(n\,c_1(X)^2-(n+2)\,c_2(X)\right)\cdot[\omega]^{n-2}.
\]

\begin{proposition}\label{11246}
	Let $(X, \omega)$ be a Hermitian manifold. The following identity holds
	\[
	\frac{n(n-1)}{n+1}\,s^2-\mathcal L_\nabla
	=\frac{n+2}{n+1}\,\mathcal D_\nabla
	+\frac{n^2}{n+1}\mid\!\rho^{\circ}_\nabla\!\mid^2.
	\]
\end{proposition}

For a K\"ahler-Einstein manifold $(X, \omega)$, the Chen-Ogiue formula \cite{CO} gives
\begin{equation}\label{11241}
	\mathcal{D}_\nabla=(n+1)\mid\!B_{\nabla}\!\mid^2,
\end{equation}
where $B_{\nabla}$ is the Bochner curvature of $\nabla$.

If $X$ is compact, $K_X$ is ample, and $\rho_\nabla=-\omega$, then $[\omega]=2\pi\,c_1(K_X)$. Integrating \eqref{11241} gives the Miyaoka-Yau inequality in its numerical form \cite{Miyaoka, Yau}
\begin{equation}\label{11247}
	\left(2(n+1)\,c_2(X)-n\,c_1(X)^2\right)\cdot c_1\left(K_X\right)^{n-2}\geq0.
\end{equation}
Equality in \eqref{11247} holds precisely when $\omega$ has negative constant holomorphic sectional curvature. Li's lower bound for the Calabi functional holds in every K\"ahler class. If the class contains a cscK metric, the bound gives an intersection inequality; in the ample canonical class it recovers \eqref{11247}.

Motivated by \eqref{11241} and \eqref{11247}, we first seek an upper bound for $\mathcal D_\nabla$ valid for every Hermitian metric and ask what equality implies about its Chern torsion. We then ask under what geometric conditions $\mathcal D_\nabla$ is nonnegative. Proposition \ref{11246} relates both questions to bounds for $\mathcal L_\nabla$.

Let $\hat R_\nabla$ denote the K\"ahler symmetrization of $R_\nabla$ and let $\nabla''T$ be the (0,1) -part of $\nabla T$.
Set
\[
\mathcal E=(n+1)\left(\mid\!\hat R_\nabla\!\mid^2+\mid\!\nabla''T\!\mid^2\right),
\]
and let $\mathcal N$ be the nonnegative function in Definition \ref{11192}.

For a function $f\in\mathcal{A}^{0}(X)$, we write $\int_X f=\int_X f\,\frac{\omega^n}{n!}$.
\begin{theorem}\label{11181}
	Let $(X,\omega)$ be a Hermitian manifold. The following hold
	\begin{equation}\label{11191}
		\mathcal D_\nabla=\mathcal E-\mathcal N,
		\quad
		\mathcal N\geq \frac{1}{2}\,(n+1)\mid\!\nabla''T\!\mid^2.
	\end{equation}
	For $n\geq 3$, the coefficient $\tfrac{1}{2}\,(n+1)$ is optimal.
	
	If $X$ is compact, then $\int_X\mathcal N=0$ if and only if
	\[
	\rho_\nabla=0, \quad \nabla T=0.
	\]
	In this case, $\omega$ is balanced, all four Chern-Ricci curvatures vanish, and $R_\nabla=B_\nabla$.
\end{theorem}

The curvature decomposition also gives
\begin{equation}\label{112610}
	\mathcal D_\nabla\leq(n+1)\mid\!R_\nabla\mid\!^2\leq\mathcal E.
\end{equation}
At $x\in X$, equality in the first inequality of \eqref{112610} holds precisely when $\rho_\nabla(x)=\kappa_\nabla(x)=0$, and equality in the second holds precisely when $R_\nabla(x)=\hat R_\nabla(x)$. Thus $\mathcal D_\nabla=\mathcal E$ at $x$ if and only if both conditions hold there. Example \ref{11248} realizes $\mathcal D_\nabla=\mathcal E$ on a compact non-K\"ahler threefold with nonzero parallel Chern
torsion and nonzero Chern curvature.

Equality in the sharp lower bound for $\mathcal N$ means
\[
\mathcal N=\frac{1}{2}(n+1)\mid\!\nabla''T\!\mid^2.
\]
On a balanced threefold, this equality holds at $x$ if and only if $\rho_\nabla(x)=\kappa_\nabla(x)=0$. If $X$ is compact and equality holds everywhere, then $\omega$ is balanced. If $X$ also admits an astheno-K\"ahler metric, equality everywhere is equivalent to $\omega$ being Ricci-flat K\"ahler.

We next consider the sign of $\mathcal D_\nabla$ for metrics satisfying the weak second Chern-Einstein equation.

\begin{theorem} \label{11224}
	Let $(X,\omega)$ be a connected Hermitian manifold satisfying the weak second Chern-Einstein equation
	\begin{equation}\label{11239}
		\kappa_\nabla=\frac{s}{n}\,\omega.
	\end{equation}
	\begin{enumerate}
		\item If $\omega$ is balanced, then
		\[
		0\leq\mathcal D_\nabla\leq\mathcal E.
		\]
		\item If $\omega$ is Gauduchon, then
		\[
		\frac{1}{n-1}\mid\!\tau\!\mid^4-\frac{2}{n}s\mid\!\tau\!\mid^2
		\leq\mathcal D_\nabla\leq\mathcal E.
		\]
		In particular, $\mathcal D_\nabla(x)\geq0$ at every point $x$ where $s(x)\leq0$.
		\item If $X$ is compact and $s$ is a nonpositive constant, then
		\[
		0\leq\int_X\mathcal D_\nabla\leq\int_X\mathcal E.
		\]
	\end{enumerate}
	At every point $x\in X$ where $\kappa_\nabla(x)=0$, one has
	\[
	\mathcal D_\nabla(x)\geq\mid\!R_\nabla\!\mid^2(x).
	\]
	In case (1), $\mathcal D_\nabla\equiv0$ if and only if $\omega$ is Chern-flat or K\"ahler with nonzero constant holomorphic sectional curvature. In case (2), if $s\leq0$ on $X$, then $\mathcal D_\nabla\equiv0$ if and only if $\omega$ is Chern-flat or K\"ahler with negative constant holomorphic sectional curvature. In case~(3), the equality
	\[
	\int_X\mathcal D_\nabla=0
	\]
	characterizes Chern-flatness when $s=0$ and K\"ahler metrics of
	negative constant holomorphic sectional curvature when $s<0$.
\end{theorem}

On compact Gauduchon surfaces the K\"ahler rigidity of a weakly second Chern-Einstein metric (following the terminology in \cite{ACS}) with nonpositive Einstein factor goes back to Gauduchon and Ivanov \cite{GI}; see also \cite{XY25}.

In Theorem \ref{11224}, the torsion enters through $s-\hat s$, where $\hat s$
is the second Chern scalar curvature. For a Gauduchon metric, $s-\hat s=\mid\!\tau\!\mid^2$.

At every point where Theorem \ref{11224} gives $\mathcal D_\nabla\geq0$,
Proposition \ref{11246} yields
\[
\mathcal L_\nabla+\frac{n^2}{n+1}\mid\!\rho_\nabla^\circ\!\mid^2
\leq \frac{n(n-1)}{n+1}s^2.
\]
The difference between the two sides is $\tfrac{n+2}{n+1}\mathcal D_\nabla$. The Gauduchon estimate in part (2) and the integral estimate in part (3) give the corresponding bounds for $\mathcal L_\nabla$ by the same identity.

We now turn from the weak second Chern-Einstein equation to balanced metrics with vanishing first Chern-Ricci form.
\begin{theorem}\label{11234}
	Let $(X,\omega)$ be a balanced manifold with $\rho_\nabla=0$. The following inequality holds
	\[
	\frac{1}{n+1}\mathcal{D}_\nabla
	\leq\mid\!B_\nabla\!\mid^2+\frac12\mid\!\nabla''T\!\mid^2.
	\]
	At $x\in X$, equality holds if and only if
	\begin{equation}\label{11235}
		\nabla_{\bar l}T_{ki\bar j}=\nabla_{\bar j}T_{ki\bar l}, \quad
		1\leq i,j,k,l\leq n,
	\end{equation}
	where $\nabla_{\bar l}T_{ki\bar j}$ are the components of $\nabla''T$. When $n=3$, condition \eqref{11235} at $x$ is equivalent to $\kappa_\nabla(x)=0$. If $X$ is compact and admits an astheno-K\"ahler metric $\chi$, then \eqref{11235} holds at every point of $X$ if and only if $\omega$ is Ricci-flat K\"ahler.
\end{theorem}

When $\rho_\nabla=0$, Proposition \ref{11246} gives
\[
\mathcal L_\nabla=-\frac{n+2}{n+1}\mathcal D_\nabla.
\]
Theorem \ref{11234} therefore gives a lower bound for $\mathcal L_\nabla$ with the same equality condition \eqref{11235}.

For Chern number identities on compact complex surfaces, see \cite{XY26}.
The next result bounds $\mathcal D_\nabla$ on a Gauduchon surface in terms of the minimum and maximum of the Chern holomorphic sectional curvature at each point.

Fix $x\in X$. Let
\[
S_x=\left\{v\in T^{1,0}_xX:\mid\!v\!\mid^{2}=1\right\}
\]
and write
\[
h_x(v)=R_{i\bar j k\bar l}(x)v^i\bar v^jv^k\bar v^l
\]
for $v\in S_x$, where $R_{i\bar{j}k\bar{l}}$ are the components of $R_{\nabla}$.

On a Hermitian surface set
\[
\mathcal H_{\min}(x)=\min_{v\in S_x}h_x(v),\quad
\mathcal H_{\max}(x)=\max_{v\in S_x}h_x(v),
\]
and
\[
\mathcal{I}_x
=\left[
3\,\mathcal H_{\min}(x)+\tfrac12\mid\!\tau\!\mid^2(x),\
3\,\mathcal H_{\max}(x)+\tfrac12\mid\!\tau\!\mid^2(x)
\right].
\]
Here $\operatorname{dist}\left(0, \mathcal{I}_x\right)
=\min_{t\in\mathcal{I}_x}\mid\!t\!\mid$.

\begin{theorem}\label{11194}
	Let \((X,\omega)\) be a connected Gauduchon surface. The following inequality holds
	\begin{equation}\label{11232}
		\mathcal{D}_\nabla+\operatorname{dist}\left(0, \mathcal{I}_x\right)^2
		\leq\mathcal{E}.
	\end{equation}
	Equality in \eqref{11232} holds at every point of $X$ if and only if $\omega$ is Ricci-flat K\"ahler or K\"ahler with nonzero constant holomorphic sectional curvature.
\end{theorem}

In complex dimension $n\geq3$, balancedness gives estimates for both densities in terms of the Bochner curvature and the first and second Chern-Ricci forms.

\begin{theorem}\label{11243}
	Let $(X,\omega)$ be a balanced manifold with $n\geq3$. The following inequalities hold
	\begin{align}
		\mathcal D_\nabla
		& \geq (n+1)\mid\!B_\nabla\!\mid^2+
		\frac{1}{n+2}\left((2n+1)\mid\!\rho^{\circ}_\nabla\!\mid^2-(n+1)^2\mid\!\kappa^{\circ}_\nabla\!\mid^2\right)\notag\\
		& \quad +\frac{2(n^2-1)}{n(n-2)(n+2)}\mid\!\kappa_\nabla-\rho_\nabla\!\mid^2 ,\label{11244}\\
		\mathcal L_\nabla
		& \leq -(n+2)\mid\!B_\nabla\!\mid^2+(n+1)\left(\mid\!\kappa^{\circ}_\nabla\!\mid^2-\mid\!\rho^{\circ}_\nabla\!\mid^2\right)+\frac{n(n-1)}{n+1}s^2-\mid\!\nabla''T\!\mid^2.\label{11245}
	\end{align}
\end{theorem}

\begin{remark}
	Let $(X,\omega)$ be a balanced manifold with $n\geq3$. If $\rho_\nabla=\kappa_\nabla$, the pointwise Li inequality
	\[
	\mathcal L_\nabla\leq\frac{n(n-1)}{n+1}s^2
	\]
	has the exact remainder
	\begin{equation}\label{11253}
		\frac{n(n-1)}{n+1}s^2-\mathcal L_\nabla
		=(n+2)\left(\mid\!B_\nabla\!\mid^2+\mid\!R_\nabla-\hat R_\nabla\!\mid^2\right).
	\end{equation}
	Equality holds if and only if
	\[
	B_\nabla=0,\quad \nabla''T=0.
	\]
\end{remark}

Under the hypotheses of Theorem \ref{11243}, assume in addition that $X$ is compact.
Integrating \eqref{11245} gives the following balanced analogue of Li's inequality. If
\[
\int_X\mid\!\kappa_\nabla^\circ\!\mid^2
\leq\int_X\mid\!\rho_\nabla^\circ\!\mid^2,
\]
then
\[
\int_X\mathcal L_\nabla
\leq\frac{n(n-1)}{n+1}\int_X s^2.
\]

Fino and Vezzoni conjectured that a compact complex manifold admitting both a balanced metric and a pluriclosed metric admits a K\"ahler metric \cite[Problem 3]{FV}. If $X$ is a compact complex threefold carrying a pluriclosed metric
$\chi$, Proposition \ref{11249} shows that a given balanced metric $\omega$ is K\"ahler if and only if
\[
\int_X(\kappa_\nabla-\rho_\nabla)\wedge\omega\wedge\chi=0.
\]

Section \ref{11072} establishes the Chern-Weil identities and their K\"ahler specializations. Section \ref{11073} proves Theorems \ref{11181}-\ref{11194}. The elliptic bundle in Example \ref{11248} also illustrates the dependence of the weighted Chern integral on the metric: its second Bott-Chern class vanishes, while the integral of its second Chern form against the balanced metric is positive. Section \ref{112410} proves the balanced estimates of Theorem \ref{11243} and the threefold criterion of Proposition \ref{11249}. It concludes with an angular integral formula for the Bochner curvature of a Hermitian metric, related to Du's formulas for compact K\"ahler-Einstein manifolds \cite{Du}.

\section{Preliminaries}\label{11072}
Let $X$ be a complex manifold with $\dim_{\mathbb{C}}X=n\geq2$. Let
\[
\omega=\sqrt{-1}g_{i\bar{j}}\mathrm{d}z^{i}\wedge\mathrm{d}\bar{z}^{j}
\]
be a Hermitian metric on $X$.

Denote the complex Laplacian of $\omega$ by $\Delta=\sqrt{-1}\Lambda_{\omega}\partial\bar{\partial}$. Denote
\[
\Gamma^{p}_{ik}=g^{p\bar{j}}\partial_{k}g_{i\bar{j}}.
\]
The Chern connection $\nabla$ is the unique Hermitian connection on $T^{1,0}X$ whose $(0,1)$-part is $\bar\partial$. Locally, we have
\[
\nabla_{k}\left(\frac{\partial}{\partial z^{i}}\right)
=\Gamma^{p}_{ik}\frac{\partial}{\partial z^{p}}, \quad
\nabla_{\bar{l}}\left(\frac{\partial}{\partial z^{i}}\right)=0,
\]
where $\nabla_{k}=\nabla_{\frac{\partial}{\partial z^{k}}}$ and $\nabla_{\bar{l}}=\nabla_{\frac{\partial}{\partial\bar{z}^{l}}}$ denote the covariant derivatives of $\nabla$ along the vector fields $\frac{\partial}{\partial z^{k}}$ and $\frac{\partial}{\partial\bar{z}^{l}}$.

For every $\theta=\theta_{\bar{l}}\mathrm{d}\bar{z}^{l}\in\mathcal{A}^{0,1}(X)$, we have
\[
\partial_{k}\theta_{\bar{l}}=\nabla_{k}\theta_{\bar{l}},\quad \partial_{\bar{l}}\bar{\theta}_{\bar{k}}=\nabla_{\bar{l}}\bar{\theta}_{\bar{k}}.
\]
These terms may be used interchangeably.

The torsion tensor of the Chern connection $\nabla$ is given by
\[
T=T_{ki\bar{j}}\mathrm{d}z^{k}\otimes\mathrm{d}z^{i}\otimes\mathrm{d}\bar{z}^{j},
\]
where
\begin{equation}\label{11033}
	T_{ki\bar{j}}=\partial_{k}g_{i\bar{j}}-\partial_{i}g_{k\bar{j}}.
\end{equation}
Denote
\[
T^{p}_{ki}=g^{p\bar{j}}T_{ki\bar{j}}=\Gamma^{p}_{ik}-\Gamma^{p}_{ki}
\]
and
\[
T_{k}=T^{i}_{ki}=g^{i\bar{j}}T_{ki\bar{j}}.
\]
Then the torsion (1,0)-form is given by
\[
\tau=\Lambda_{\omega}\partial\omega=T_{k}\mathrm{d}z^{k}.
\]

Since each Hermitian metric $\omega$ satisfies
\begin{equation}\label{06041}
	\partial\omega^{n-1}=\tau\wedge\omega^{n-1},
\end{equation}
$\omega$ is Gauduchon if and only if (e.g., section 7.1 in \cite{LT})
\begin{equation}\label{10111}
	0=\mid\!\tau\!\mid^{2}+g^{k\bar{l}}\nabla_{\bar{l}}T_{k}.
\end{equation}
where
\[
\mid\!\tau\!\mid^{2}=g^{k\bar{l}}T_{k}\bar{T}_{l}, \quad g^{k\bar{l}}\nabla_{\bar{l}}T_{k}=g^{k\bar{l}}\nabla_{k}\bar{T}_{l}.
\]
And $\omega$ is balanced if $\tau=0$.

In particular, for $n=2$, we obtain from \eqref{11033} and \eqref{06041} that
\begin{equation}\label{11034}
	T_{ki\bar{j}}=T_{k}g_{i\bar{j}}-T_{i}g_{k\bar{j}}.
\end{equation}

Integration by parts will be used. Let $(X, \omega)$ be a compact Hermitian manifold. For $\theta=\theta_{\bar{l}}\mathrm{d}\bar{z}^{l}\in\mathcal{A}^{0,1}(X)$, we obtain from \eqref{06041} that
\[
\begin{split}
	0
	& = n\int_{X}\sqrt{-1}\mathrm{d}\left(\left(\theta-\bar{\theta}\right)\wedge\omega^{n-1}\right)\\
	& = n\int_{X}\sqrt{-1}\left(\partial\theta-\bar{\partial}\bar{\theta}+\tau\wedge\theta+\bar{\theta}\wedge\bar{\tau}\right)\wedge\omega^{n-1}\\
	& = \int_{X}2\operatorname{Re}\left(g^{k\bar{l}}\left(\nabla_{k}\theta_{\bar{l}}+T_{k}\theta_{\bar{l}}\right)\right)\omega^{n}.
\end{split}
\]

The curvature $R_{\nabla}$ of the Chern connection $\nabla$ of a Hermitian metric $\omega$ is
\[
R_{\nabla}=\nabla\circ\nabla\in\mathcal{A}^{1,1}\left(\mathrm{End}\left(T^{1,0}X\right)\right).
\]
Locally, we calculate
\[
R_{\nabla}=\bar{\partial}\left(\partial gg^{-1}\right)^{T}
=R^{p}_{ik\bar{l}}\,\mathrm{d}z^{k}\wedge\mathrm{d}\bar{z}^{l}\otimes\mathrm{d}z^{i}\otimes\frac{\partial}{\partial z^{p}},
\]
where (see page 12 in \cite{Kobayashi})
\[
R^{p}_{ik\bar{l}}=g^{p\bar{j}}\left(-\partial_{\bar{l}}\partial_{k}g_{i\bar{j}}+g^{m\bar{n}}\partial_{\bar{l}}g_{m\bar{j}}\partial_{k}g_{i\bar{n}}\right).
\]
Denote
\begin{equation}\label{09221}
	R_{i\bar{j}k\bar{l}}=g_{p\bar{j}}R^{p}_{ik\bar{l}}
	=-\partial_{\bar{l}}\partial_{k}g_{i\bar{j}}+g^{p\bar{q}}\partial_{\bar{l}}g_{p\bar{j}}\partial_{k}g_{i\bar{q}}.
\end{equation}
Then we have $R_{i\bar{j}k\bar{l}}=\bar{R}_{j\bar{i}l\bar{k}}$. The following identities are well-known (e.g., (2.6) in \cite{SW}):
\begin{equation}\label{11022}
	R_{i\bar{j}k\bar{l}}-R_{k\bar{j}i\bar{l}}
	=\nabla_{\bar{l}}T_{ik\bar{j}}, \quad
	R_{i\bar{j}k\bar{l}}-R_{k\bar{l}i\bar{j}}
	=\nabla_{\bar{j}}T_{ik\bar{l}}+\nabla_{k}\bar{T}_{jl\bar{i}}.
\end{equation}

There are four Chern-Ricci curvatures of $\nabla$:
\[
R_{i\bar{j}}=g^{k\bar{l}}R_{k\bar{l}i\bar{j}}, \quad 
K_{i\bar{j}}=g^{k\bar{l}}R_{i\bar{j}k\bar{l}}, \quad
R^{3}_{i\bar{j}}=g^{k\bar{l}}R_{k\bar{j}i\bar{l}}, \quad
R^{4}_{i\bar{j}}=g^{k\bar{l}}R_{i\bar{l}k\bar{j}}.
\]
There are two scalar curvatures:
\[
s=g^{i\bar{j}}R_{i\bar{j}}=g^{i\bar{j}}K_{i\bar{j}}, \quad
\hat{s}=g^{i\bar{j}}R^{3}_{i\bar{j}}=g^{i\bar{j}}R^{4}_{i\bar{j}}.
\]
In particular, if $\omega$ is Gauduchon, we obtain from \eqref{10111} and \eqref{11022} that
\begin{equation}\label{10271}
	s-\hat{s}=\mid\!\tau\!\mid^{2}.
\end{equation}

Denote the two real curvature (1,1)-forms:
\[
\rho_{\nabla}=\sqrt{-1}\, R_{i\bar{j}}\mathrm{d}z^{i}\wedge\mathrm{d}\bar{z}^{j}, \quad
\kappa_{\nabla}=\sqrt{-1}\, K_{i\bar{j}}\mathrm{d}z^{i}\wedge\mathrm{d}\bar{z}^{j}.
\]
Then
\[
\rho^{\circ}_{\nabla}
=\sqrt{-1}\, R^{\circ}_{i\bar{j}}\mathrm{d}z^{i}\wedge\mathrm{d}\bar{z}^{j}, \quad
\kappa^{\circ}_{\nabla}
=\sqrt{-1}\, K^{\circ}_{i\bar{j}}\mathrm{d}z^{i}\wedge\mathrm{d}\bar{z}^{j},
\]
where
\[
R^\circ_{i\bar{j}}=R_{i\bar{j}}-\frac{s}{n}g_{i\bar{j}}, \quad
K^\circ_{i\bar{j}}=K_{i\bar{j}}-\frac{s}{n}g_{i\bar{j}}
\]

The K\"{a}hler symmetrization $\hat{R}_{\nabla}$ of $R_{\nabla}$ was defined in \cite{B}, whose components are
\begin{equation}\label{11211}
	\hat{R}_{i\bar{j}k\bar{l}}=\frac{1}{4}\left(R_{i\bar{j}k\bar{l}}+R_{k\bar{l}i\bar{j}}+R_{k\bar{j}i\bar{l}}+R_{i\bar{l}k\bar{j}}\right).
\end{equation}
While its doubly alternating components, obtained by alternating separately in the holomorphic and antiholomorphic indices, are
\[
\tilde{R}_{i\bar j k\bar l}
=\frac14\left(R_{i\bar j k\bar l}+R_{k\bar l i\bar j}-R_{k\bar j i\bar l}-R_{i\bar l k\bar j}\right).
\]
Denote
\begin{equation}\label{11197}
	M_{i\bar j k\bar l}=\frac{1}{2}\left(R_{i\bar{j}k\bar{l}}-R_{k\bar{l}i\bar{j}}\right).
\end{equation}
Direct clculation yields the following.
\begin{proposition}\label{11221}
	Let $(X, \omega)$ be a Hermitian manifold. There is a orthogonal decomposition
	\[
	R_{i\bar j k\bar l}=\hat{R}_{i\bar{j}k\bar{l}}+\tilde{R}_{i\bar{j}k\bar{l}}+M_{i\bar j k\bar l}.
	\]
\end{proposition}

Denote
\begin{equation}\label{07101}
	\hat{R}_{i\bar{j}}=g^{k\bar{l}}\hat{R}_{k\bar{l}i\bar{j}}=\frac{1}{4}(R_{i\bar{j}}+K_{i\bar{j}}+R^{3}_{i\bar{j}}+R^{4}_{i\bar{j}}).
\end{equation}
The trace-free part is
\begin{equation*}
	\hat{R}_{i\bar{j}}^{\circ}=\hat{R}_{i\bar{j}}-\frac{s+\hat{s}}{2n}g_{i\bar{j}}.
\end{equation*}
\begin{definition}\label{09131}
	Let $(X, \omega)$ be a Hermitian manifold. The Bochner curvature $B_{\nabla}$ of $\nabla$ is
	\[
	\begin{split}
		B_{i\bar{j}k\bar{l}}
		&=\hat{R}_{i\bar{j}k\bar{l}}- \frac{1}{n+2}\left(\hat{R}_{i\bar{j}}^{\circ}g_{k\bar{l}}+\hat{R}_{i\bar{l}}^{\circ}g_{k\bar{j}}+\hat{R}_{k\bar{j}}^{\circ}g_{i\bar{l}}+\hat{R}_{k\bar{l}}^{\circ}g_{i\bar{j}}\right)\\
		& \quad -\frac{s+\hat{s}}{2n(n+1)}\left(g_{i\bar{j}}g_{k\bar{l}}+g_{i\bar{l}}g_{k\bar{j}}\right).
	\end{split}
	\]
\end{definition}

For $1\leq p,q\leq n$, the norm of a tensor
\[
\mathcal{T}=\mathcal{T}_{i_{1}\cdots i_{p}\bar{j}_{1}\cdots \bar{j}_{q}}\mathrm{d}z^{i_{1}}\otimes\cdots\otimes\mathrm{d}z^{i_{p}}\otimes\mathrm{d}\bar{z}^{j_{1}}\otimes\cdots\otimes\mathrm{d}\bar{z}^{j_{q}}
\]
is denoted by
\[
\mid\!\mathcal{T}_{i_{1}\cdots i_{p}\bar{j}_{1}\cdots \bar{j}_{q}}\!\mid^{2}
=g^{i_{1}\bar{l}_{1}}\cdots g^{i_{p}\bar{l}_{p}}g^{k_{1}\bar{j}_{1}}\cdots g^{k_{q}\bar{j}_{q}}\mathcal{T}_{i_{1}\cdots i_{p}\bar{j}_{1}\cdots \bar{j}_{q}}
\bar{\mathcal{T}}_{l_{1}\cdots l_{p}\bar{k}_{1}\cdots \bar{k}_{q}}.
\]

\begin{proof}[Proof of Proposition \ref{11246}]
	Let $(X, \omega)$ be a Hermitian manifold. By \eqref{11233} and the calculations in \cite{Kobayashi} (Chap. 4.4), we have
	\begin{equation}\label{11222}
		\begin{split}
			\mathcal{D}_\nabla
			&=(n+1)\left(\mid\!R_{i\bar{j}k\bar{l}}\!\mid^2-\mid\!K_{i\bar{j}}\!\mid^2\right)+s^{2}-\mid\!R_{i\bar{j}}\!\mid^2,\\
			\mathcal{L}_\nabla
			&=(n+2)\left(\mid\!K_{i\bar{j}}\!\mid^2-\mid\!R_{i\bar{j}k\bar{l}}\!\mid^2\right)+(n-2)\left(s^2-\mid\!R_{i\bar{j}}\!\mid^2\right).
		\end{split}
	\end{equation}
	Canceling the term $\mid\!R_{i\bar{j}k\bar{l}}\!\mid^2$ yields the conclusion.
\end{proof}

When $\omega$ is K\"ahler, by Proposition \ref{11246} and Definition \ref{09131},
\[
\begin{split}
	\mathcal{D}_\nabla
	&=(n+1)\mid\!B_{i\bar{j}k\bar{l}}\!\mid^2-\frac{n^2}{n+2}\mid\!R^\circ_{i\bar{j}}\!\mid^2,\\
	\mathcal{L}_\nabla
	&=\frac{n(n-1)}{n+1}\,s^2-(n+2)\mid\!B_{i\bar{j}k\bar{l}}\!\mid^2.
\end{split}
\]
The second equality gives Li's bound for the Calabi functional upon integration \cite{Li}.

Suppose $(X, \omega)$ is compact cscK. We have
\[
\begin{split}
	s&=2\pi\,n\frac{c_1(X)\cdot[\omega]^{n-1}}{[\omega]^n},\\
	\int_X\mid\!R^\circ_{i\bar{j}}\!\mid^2
	&=
	\frac{4\pi^2}{(n-2)!}
	\left(\frac{\left(c_1(X)\cdot[\omega]^{n-1}\right)^2}{[\omega]^n}-c_1(X)^2\cdot[\omega]^{n-2}\right).
\end{split}
\]
Hence Li's intersection inequality may be written as
\[
\int_X\mathcal D_\nabla
\geq
-\frac{\left(2\pi\,n\right)^2}{(n+2)(n-2)!}
\left(\frac{\left(c_1(X)\cdot[\omega]^{n-1}\right)^2}{[\omega]^n}-c_1(X)^2\cdot [\omega]^{n-2}\right).
\]
The term in parentheses vanishes when $c_1(X)$ is proportional to $[\omega]$. If $K_X$ is ample, the K\"ahler-Einstein metric in the class $[\omega]=c_1\left(K_X\right)$ satisfies $\rho_\nabla=-2\pi\omega$. The preceding inequality then gives \eqref{11247}.

In complex dimension two, the integrals of $\mathcal D_\nabla$ and $\mathcal L_\nabla$ are linear combinations of Chern numbers. In higher dimensions they use the weight $\frac{\omega^{n-2}}{(n-2)!}$, which need not be $\partial\bar\partial$-closed even when $\omega$ is balanced (see Example \ref{11248}).

We end this section by the definitions of the nonnegative term $\mathcal{N}$ in Theorem \ref{11181}.
\begin{definition}\label{11192}
	Let $(X, \omega)$ be a Hermitian manifold. Define
	\[
	\mathcal{N}
	=(n+1)\left(3\mid\!\tilde{R}_{i\bar{j}k\bar{l}}\!\mid^2+\mid\!M_{i\bar{j}k\bar{l}}\!\mid^2+\mid\!K^\circ_{i\bar{j}}\!\mid^2\right)+\mid\!R^\circ_{i\bar{j}}\!\mid^2+\frac{2}{n}s^2.
	\]
\end{definition}

\section{Proofs of Theorems \ref{11181} to \ref{11194}}\label{11073}

In this section, we first express $\mid\!\nabla''T\!\mid^2$ in terms of $R_\nabla$ and use this identity to prove the sharp lower bound for $\mathcal N$ in Theorem \ref{11181}. We then establish estimates for $\mathcal D_\nabla$ and characterize equality for weakly second Chern-Einstein metrics, for balanced metrics with $\rho_\nabla=0$, and for Gauduchon surfaces.
\begin{lemma}\label{11195}
	Let $(X,\omega)$ be a Hermitian manifold. The following identity holds
	\[
	\mid\!\nabla_{\bar{l}}T_{ki\bar{j}}\!\mid^2
	=4\mid\!\tilde{R}_{i\bar{j}k\bar{l}}\!\mid^2+2\mid\!M_{i\bar{j}k\bar{l}}\!\mid^2.
	\]
\end{lemma}
\begin{proof}
	By \eqref{11022} and \eqref{11197}, we have
	\begin{align*}
		\mid\!\tilde{R}_{i\bar{j}k\bar{l}}\!\mid^2
		&=\frac{1}{16}\left(2\mid\!\nabla_{\bar{l}}T_{ki\bar{j}}\!\mid^{2}
		+2\operatorname{Re}\left(g^{i\bar{q}}g^{p\bar{j}}g^{k\bar{n}}g^{m\bar{l}}\nabla_{\bar{l}}T_{ki\bar{j}}\nabla_{p}\bar{T}_{qn\bar{m}}\right)\right),\\
		\mid\!M_{i\bar{j}k\bar{l}}\!\mid^2
		&=\frac{1}{4}\left(2\mid\!\nabla_{\bar{l}}T_{ki\bar{j}}\!\mid^{2}
		-2\operatorname{Re}\left(g^{i\bar{q}}g^{p\bar{j}}g^{k\bar{n}}g^{m\bar{l}}\nabla_{\bar{l}}T_{ki\bar{j}}\nabla_{\bar{q}}T_{pm\bar{n}}\right)\right).
	\end{align*}
	Then we obtain
	\[
	\begin{split}
		4\mid\!\tilde{R}_{i\bar{j}k\bar{l}}\!\mid^2+2\mid\!M_{i\bar{j}k\bar{l}}\!\mid^2
		& = \frac{3}{2}\mid\!\nabla_{\bar{l}}T_{ki\bar{j}}\!\mid^2
		+\frac{1}{4}\,2\operatorname{Re}\left(g^{i\bar{q}}g^{p\bar{j}}g^{k\bar{n}}g^{m\bar{l}}\nabla_{\bar{l}}T_{ki\bar{j}}\nabla_{p}\bar{T}_{qn\bar{m}}\right)\\
		& \quad -\frac{1}{2}\,2\operatorname{Re}\left(g^{i\bar{q}}g^{p\bar{j}}g^{k\bar{n}}g^{m\bar{l}}\nabla_{\bar{l}}T_{ki\bar{j}}\nabla_{\bar{q}}T_{pm\bar{n}}\right)\\
		& = \mid\!\nabla_{\bar{l}}T_{ki\bar{j}}\!\mid^2
		+\frac{1}{4}\,2\operatorname{Re}\left(g^{i\bar{q}}g^{p\bar{j}}g^{k\bar{n}}g^{m\bar{l}}\nabla_{\bar{l}}T_{ki\bar{j}}\left(\nabla_{\bar{n}}T_{mp\bar{q}}-\nabla_{\bar{q}}T_{pm\bar{n}}\right)\right),
	\end{split}
	\]
	where we have used \eqref{11022} for the second equality. While it is direct to check (e.g., under a unitary frame) that
	\[
	2\operatorname{Re}\left(g^{i\bar{q}}g^{p\bar{j}}g^{k\bar{n}}g^{m\bar{l}}\nabla_{\bar{l}}T_{ki\bar{j}}\left(\nabla_{\bar{n}}T_{mp\bar{q}}-\nabla_{\bar{q}}T_{pm\bar{n}}\right)\right)=0
	\]
	The desired identity follows.
\end{proof}

\begin{proof}[Proof of Theorem \ref{11181}]
	Let $(X, \omega)$ be a Hermitian manifold. By Proposition \ref{11221},
	\begin{equation}\label{11227}
		\mid\!R_{i\bar{j}k\bar{l}}\!\mid^2=\mid\!\hat{R}_{i\bar{j}k\bar{l}}\!\mid^2+\mid\!\tilde{R}_{i\bar{j}k\bar{l}}\!\mid^2+\mid\!M_{i\bar{j}k\bar{l}}\!\mid^2.
	\end{equation}
	By Definition \ref{11192}, Lemma \ref{11195} and \eqref{11227},
	\[
	\begin{split}
		\mathcal{E}-\mathcal{N}
		& = (n+1)\left(\mid\!R_{i\bar{j}k\bar{l}}\!\mid^2-\mid\!K_{i\bar{j}}^{\circ}\!\mid^2\right)-\mid\!R_{i\bar{j}}^{\circ}\!\mid^2-\frac{2}{n}s^{2}\\
		& = (n+1)\left(\mid\!R_{i\bar{j}k\bar{l}}\!\mid^2-\mid\!K_{i\bar{j}}\!\mid^2\right)-\mid\!R_{i\bar{j}}\!\mid^2+s^{2}.
	\end{split}
	\]
	Together with \eqref{11222}, this gives
	\[
	\mathcal D_\nabla=\mathcal E-\mathcal N.
	\]
	To prove the lower bound for $\mathcal N$, we use Definition \ref{11192} and Lemma \ref{11195} to obtain
	\[
	\begin{aligned}
		\mathcal N-\frac{1}{2}\,(n+1)\mid\!\nabla''T\!\mid^2
		&=(n+1)\mid\!\tilde R_{i\bar j k\bar l}\!\mid^2
		+(n+1)\mid\!K^\circ_{i\bar j}\!\mid^2
		+\mid\!R^\circ_{i\bar j}\!\mid^2+\frac{2}{n}s^2\\
		&\geq 0.
	\end{aligned}
	\]
	Equality in this lower bound holds at $x\in X$ if and only if
	\[
	\tilde R_{i\bar j k\bar l}(x)=0,\quad
	R_{i\bar j}(x)=K_{i\bar j}(x)=0.
	\]
	By \eqref{11022},
	\[
	\tilde R_{i\bar j k\bar l}
	=\frac{1}{4}\left(\nabla_{\bar l}T_{ik\bar j}-\nabla_{\bar j}T_{ik\bar l}\right).
	\]
	Thus $\tilde R_{i\bar j k\bar l}(x)=0$ if and only if
	\[
	\nabla_{\bar l}T_{ik\bar j}=\nabla_{\bar j}T_{ik\bar l},\quad
	1\leq i,j,k,l\leq n.
	\]
	
	To prove optimality, let $n\geq3$. In a neighborhood of $0\in\mathbb C^n$, pick the metric
	\[
	\left(g_{i\bar j}\right)
	=
	\begin{pmatrix}
		1 & 0 & -z^2\bar{z}^3 & 0\\
		0 & 1 &  z^1\bar{z}^3 & 0\\
		-z^3\bar{z}^2 & z^3\bar{z}^1 & 1 & 0\\
		0 & 0 & 0 & I_{n-3}
	\end{pmatrix},
	\]
	whose curvature at 0 is
	\[
	R_\nabla^{T}(0)=
	\begin{pmatrix}
		0 & 0 & \mathrm dz^2\wedge\mathrm d\bar z^3 & 0\\
		0 & 0 & -\mathrm dz^1\wedge\mathrm d\bar z^3 & 0\\
		\mathrm dz^3\wedge\mathrm d\bar z^2
		&-\mathrm dz^3\wedge\mathrm d\bar z^1 & 0 & 0\\
		0 & 0 & 0 & 0_{n-3}
	\end{pmatrix}.
	\]
	Definition \ref{11192} and Lemma \ref{11195}  give
	\[
	\mathcal N=4(n+1),\quad \mid\!\nabla''T\!\mid^2=8.
	\]
	It follows that $\tfrac{1}{2}\,(n+1)$ is the largest constant for which the pointwise lower bound for $\mathcal N$.
	
	Suppose $X$ is compact and $\int_{X}\mathcal{N}=0$. Since the terms in $\mathcal{N}$ are all nonnegative, their integral vanishes if and only if
	\begin{equation}\label{11223}
		R_{i\bar{j}k\bar{l}}=\hat{R}_{i\bar{j}k\bar{l}}, \quad \tilde{R}_{i\bar{j}k\bar{l}}=0,\quad R_{i\bar{j}}=K_{i\bar{j}}=0.
	\end{equation}
	By Lemma \ref{11195}, we have $\nabla_{\bar{l}}T_{ki\bar{j}}=0$. The Bochner formula reads (using $K_{i\bar{j}}=0$ in \eqref{11223})
	\begin{equation}\label{11229}
		\begin{split}
			\Delta\mid\!T_{ki\bar{j}}\!\mid^{2}
			& = g^{i\bar{j}}g^{k\bar{l}}g^{m\bar{n}}g^{p\bar{q}}T_{ki\bar{q}}[\nabla_{m},\nabla_{\bar{n}}]\bar{T}_{lj\bar{p}}+\mid\!\nabla_{p}T_{ki\bar{j}}\!\mid^{2}\\
			& = \mid\!\nabla_{p}T_{ki\bar{j}}\!\mid^{2}\geq0.
		\end{split}
	\end{equation}
	The maximum principle gives $\nabla'T=0$ on each connected component of $X$. Hence $\nabla T=0$. Conversely, if $\rho_\nabla=0$ and $\nabla T=0$, then \eqref{11022} gives the K\"ahler symmetries of $R_\nabla$. Thus $\kappa_\nabla=\rho_\nabla=0$, and Definition \ref{11192} gives $\mathcal N=0$.
	
	Under these equivalent conditions, the torsion trace satisfies
	$\nabla_{\bar l}T_k=0$. Since $X$ is compact, integration by parts gives
	\begin{equation}\label{112310}
		\int_{X}\mid\!\tau\!\mid^{2}=-\int_{X}g^{k\bar{l}}\nabla_{\bar{l}}T_{k},
	\end{equation}
	we have $\tau=0$, i.e., $\omega$ is balanced. By \eqref{11022},
	\begin{equation}\label{11226}
		0=R_{i\bar{j}}=R_{i\bar{j}}^{3}=R_{i\bar{j}}^{4}.
	\end{equation}
	By Definition \ref{09131}, \eqref{07101} and \eqref{11223},
	\[
	B_{i\bar{j}k\bar{l}}=R_{i\bar{j}k\bar{l}}.
	\]
	The proof is complete.
\end{proof}

\begin{remark}
	By Definition \ref{11192} and \eqref{11191}, for every $c\in\mathbb R$ one has
	\[
	\begin{split}
		(n+1)\left(\mid\!\hat R_\nabla\!\mid^2+c\,\mid\!\nabla''T\!\mid^2\right)
		-\mathcal D_\nabla
		& = \left(\mathcal N-\tfrac{1}{2}\,(n+1)\mid\!\nabla''T\!\mid^2\right)\\
		& \quad +(n+1)\left(c-\tfrac12\right)\mid\!\nabla''T\!\mid^2.
	\end{split}
	\]
	For $c>\tfrac12$, the right hand side vanishes at $x$ if and only if
	\[
	\rho_\nabla(x)=\kappa_\nabla(x)=0, \quad \nabla''T(x)=0.
	\]
	By \eqref{11222}, \eqref{11227} and Lemma \ref{11195}, we get \eqref{112610}. Equality in the first inequality holds precisely where $\rho_\nabla=\kappa_\nabla=0$. Equality in the second holds precisely
	where $R_\nabla=\hat R_\nabla$.
\end{remark}

To prove Theorem \ref{11224}, we need a additional result. Let $(X, \omega)$ be a Hermitian manifold. Define a tensor with components
\begin{equation}\label{112311}
	Q_{i\bar{j}k\bar{l}}
	=R_{i\bar{j}k\bar{l}}-\frac{1}{n}\left(R_{k\bar{l}}g_{i\bar{j}}+K^{\circ}_{i\bar{j}}g_{k\bar{l}}\right)
	-\frac{1}{n^{2}-1}\left(\hat{s}-\frac{s}{n}\right)\left(g_{i\bar{l}}g_{k\bar{j}}-\frac{1}{n}g_{i\bar{j}}g_{k\bar{l}}\right).
\end{equation}
\begin{lemma}\label{11238}
	Let $(X, \omega)$ be a Hermitian manifold. The following identity holds
	\[
	\mathcal{D}_{\nabla}+\frac{n^{2}-1}{n}\mid\!K^{\circ}_{i\bar{j}}\!\mid^{2}
	=(n+1)\mid\!Q_{i\bar{j}k\bar{l}}\!\mid^{2}
	+\frac{1}{n}\mid\!R^{\circ}_{i\bar{j}}\!\mid^{2}
	+\frac{1}{n-1}\left(s-\hat{s}\right)^{2}-\frac{2}{n}\,s\left(s-\hat{s}\right).
	\]
\end{lemma}
\begin{proof}
	Direct calculation yields
	\[
	\begin{split}
		\mid\!Q_{i\bar{j}k\bar{l}}\!\mid^{2}
		& = \mid\!R_{i\bar{j}k\bar{l}}\!\mid^{2}-\frac{1}{n}\left(\mid\!R_{i\bar{j}}\!\mid^{2}+\mid\!K^{\circ}_{i\bar{j}}\!\mid^{2}\right)-\frac{1}{n^{2}-1}\left(\hat{s}-\frac{s}{n}\right)^{2}\\
		& = \mid\!R_{i\bar{j}k\bar{l}}\!\mid^{2}-\frac{1}{n}\left(\mid\!R^{\circ}_{i\bar{j}}\!\mid^{2}+\mid\!K^{\circ}_{i\bar{j}}\!\mid^{2}\right)-\frac{1}{n^{2}}s^{2}-\frac{1}{n^{2}-1}\left(\hat{s}-\frac{s}{n}\right)^{2}.
	\end{split}
	\]
	Combining this with \eqref{11222} yields
	\[
	\mathcal{D}_{\nabla}-(n+1)\mid\!Q_{i\bar{j}k\bar{l}}\!\mid^{2}
	=-\frac{n^{2}-1}{n}\mid\!K^{\circ}_{i\bar{j}}\!\mid^{2}+\frac{1}{n}\mid\!R^{\circ}_{i\bar{j}}\!\mid^{2}
	+\frac{1}{n-1}\left(\hat{s}-\frac{s}{n}\right)^{2}-\frac{n-1}{n^{2}}s^{2}.
	\]
	Rearranging this identity yields the conclusion.
\end{proof}
\begin{proof}[Proof of Theorem \ref{11224}]
	Let $(X, \omega)$ be a Hermitian manifold satisfying the equation \eqref{11239}.
	
	If $\omega$ ia balanced, we obtain from \eqref{11022} that $s=\hat{s}$. By Lemma \ref{11238},
	\[
	\mathcal{D}_{\nabla}
	=(n+1)\mid\!Q_{i\bar{j}k\bar{l}}\!\mid^{2}
	+\frac{1}{n}\mid\!R^{\circ}_{i\bar{j}}\!\mid^{2}\geq0.
	\]
	
	If $\omega$ is Gauduchon, \eqref{10271} and Lemma \ref{11238} give
	\[
	\mathcal D_\nabla-\frac{1}{n-1}\mid\!\tau\!\mid^4+\frac{2}{n}\, s\mid\!\tau\!\mid^2
	=(n+1)\mid\!Q_{i\bar j k\bar l}\!\mid^2+\frac1n\mid\!R^\circ_{i\bar j}\!\mid^2.
	\]
	This proves the lower bound in (2). Equality in this bound at $x$ is equivalent to
	\[
	Q_{i\bar j k\bar l}(x)=0, \quad R^\circ_{i\bar j}(x)=0.
	\]
	For $s\leq0$ the lower bound is nonnegative.
	
	Suppose now that $X$ is compact. By \eqref{11022} and \eqref{112310},
	\begin{equation}\label{11261}
		\int_X\left(s-\hat s\right)=\int_X\mid\!\tau\!\mid^2.
	\end{equation}
	If $s$ is a constant, we obtain from Lemma \ref{11238} that
	\[
	\int_X\mathcal D_\nabla
	=
	\int_X\left(
	(n+1)\mid\!Q_{i\bar j k\bar l}\!\mid^2
	+\frac1n\mid\!R^\circ_{i\bar j}\!\mid^2
	+\frac{1}{n-1}\left(s-\hat s\right)^2
	\right)-\frac{2}{n}\, s\int_X\mid\!\tau\!\mid^2.
	\]
	For $s\leq0$, the right-hand side is nonnegative. The upper bound in part (3) follows by integrating $\mathcal D_\nabla\leq\mathcal E$.
	
	At any point where $\kappa_\nabla=0$, we have $s=0$. By \eqref{11222},
	\[
	\mathcal D_\nabla-\mid\!R_\nabla\!\mid^2
	=n\left|R_{i\bar j k\bar l}
	-\frac1n g_{i\bar j}R_{k\bar l}\right|^2
	\geq0.
	\]
	In particular, at such a point, $\mathcal D_\nabla(x)=0$ if and only if $R_\nabla(x)=0$.
	
	In part (3), if $s=0$, then \eqref{11239} gives $\kappa_\nabla=0$ on $X$. It follows that $\int_X\mathcal D_\nabla=0$ if and only if $R_\nabla\equiv0$.
	
	Suppose that $\mathcal D_\nabla\equiv0$ in part (1), or that $\mathcal D_\nabla\equiv0$ and $s\leq0$ on $X$ in part (2). In part (3), suppose that $s<0$ and $\int_X\mathcal D_\nabla=0$. The corresponding identities give on $X$:
	\[
	Q_{i\bar j k\bar l}=0,\quad R^\circ_{i\bar j}=0.
	\]
	In part (1), $\tau=0$ by balancedness. In part (2), the term $\tfrac{1}{n-1}\mid\!\tau\!\mid^4$ vanishes at every point. In part (3), since $-\tfrac{2}{n}s>0$, the integral identity gives $\int_X\mid\!\tau\!\mid^2=0$. Thus in all three cases $\tau=0$ and $s=\hat s$. Substitution in \eqref{112311} gives
	\begin{equation}\label{112312}
		R_{i\bar j k\bar l}=\frac{s}{n(n+1)}\left(g_{i\bar j}g_{k\bar l}+g_{i\bar l}g_{k\bar j}\right),
	\end{equation}
	Since $\rho_\nabla=\frac{1}{n}s\,\omega$ is closed and $\mathrm d\omega^{n-1}=0$, it follows that
	\[
	0=\mathrm d\left(s\,\omega\right)\wedge\omega^{n-2}
	=\mathrm ds\wedge\omega^{n-1}.
	\]
	Wedge multiplication by $\omega^{n-1}$ is injective on one-forms, so $s$ is constant on the connected manifold $X$. If $s=0$, \eqref{112312} gives $R_\nabla=0$. If $s\neq0$, then $\mathrm d(s\omega)=0$ gives $\mathrm d\omega=0$, and \eqref{112312} gives constant holomorphic sectional curvature $\tfrac{2s}{n(n+1)}$. Conversely, either curvature tensor gives $\mathcal D_\nabla=0$ by \eqref{11222}.
\end{proof}

For a K\"ahler Einstein metric, we have
\[
\tilde{R}_{i\bar j k\bar l}=M_{i\bar j k\bar l}=0, \quad
R^\circ_{i\bar j}=K^\circ_{i\bar j}=0,\quad s=\hat s,
\]
which together with Definition \ref{09131} imply \eqref{11241}.

Every Ricci-flat K\"ahler metric satisfies $\mathcal N=0$ and hence
$\mathcal D_\nabla=\mathcal E$. For the Fubini-Study metric on $\mathbb{CP}^n$, the holomorphic sectional curvature is a positive constant, so $B_\nabla=0$ and \eqref{11241} gives $\mathcal D_\nabla=0$. Thus the Fubini-Study metric attains the lower bound in Theorem \ref{11224} (1), while
\[
\nabla''T=0,\quad
\mathcal N=\frac{2}{n}s^2>0.
\]

Without the second Chern-Einstein hypothesis, the discriminant can be
negative even for a K\"ahler surface. On the product of the projective
line and a curve of genus $g>1$, the splitting of the tangent bundle
gives
\[
\int_Xc_2(X)=4(1-g), \quad \int_Xc_1(X)^2=8(1-g).
\]

\begin{example}
	For $n\geq2$ and $(\delta_{1}, \cdots, \delta_{n})\in\mathbb{C}^{n}\setminus\{0\}$ with $\mid\!\delta_{1}\!\mid=\cdots=\mid\!\delta_{n}\!\mid\neq1$, let $S^{2n-1}\times S^{1}=(\mathbb{C}^{n}\setminus\{0\})/\sim$ be the Hopf manifold, where
	\[
	(z^{1}, \cdots, z^{n})\sim(\delta_{1}z^{1}, \cdots, \delta_{n}z^{n}).
	\]
	There is a natural Gauduchon metric
	\[
	\omega=\sqrt{-1}\,\frac{\delta_{ij}}{\mid\!z\!\mid^{2}}\,\mathrm{d}z^{i}\wedge\mathrm{d}\bar{z}^{j}.
	\]
	It is showed in \cite{Yang21} that
	\[
	\kappa_{\nabla}=(n-1)\,\omega, \quad
	\mathcal{D}_{\nabla}=-n(n-1)(n-2).
	\]
	Thus, for $n\geq3$, we have $\mathcal{D}_{\nabla}<0$.
\end{example}

\begin{proof}[Proof of Theorem \ref{11234}]
	Let $(X, \omega)$ be a balanced manifold with $\rho_{\nabla}=0$. By \eqref{11222},
	\begin{equation}\label{11228}
		\frac{1}{n+1}\mathcal{D}_{\nabla}
		=\mid\!R_{i\bar{j}k\bar{l}}\!\mid^2-\mid\!K_{i\bar{j}}\!\mid^2.
	\end{equation}
	By \eqref{07101} and \eqref{11226},
	\[
	\hat{R}_{i\bar{j}}^{\circ}=\hat{R}_{i\bar{j}}=\frac{1}{4}K_{i\bar{j}}.
	\]
	By Corollary 1 in \cite{Yang},
	\[
	\mid\!\hat{R}_{i\bar{j}k\bar{l}}\!\mid^{2}
	=\frac{1}{4(n+2)}\mid\!K_{i\bar{j}}\!\mid^{2}
	+\mid\!B_{i\bar{j}k\bar{l}}\!\mid^{2}.
	\]
	
	Inserting \eqref{11227} and the above identity into \eqref{11228} yields
	\begin{equation}\label{11213}
		\begin{split}
			\frac{1}{n+1}\,\mathcal{D}_{\nabla}
			& = \mid\!B_{i\bar{j}k\bar{l}}\!\mid^{2}+\mid\!\tilde{R}_{i\bar{j}k\bar{l}}\!\mid^2+\mid\!M_{i\bar{j}k\bar{l}}\!\mid^2-\frac{4n+7}{4(n+2)}\mid\!K_{i\bar{j}}\!\mid^{2}\\
			& = \mid\!B_{i\bar{j}k\bar{l}}\!\mid^{2}
			+\frac{1}{2}\mid\!\nabla_{\bar{l}}T_{ki\bar{j}}\!\mid^2
			-\mid\!\tilde{R}_{i\bar{j}k\bar{l}}\!\mid^2
			-\frac{4n+7}{4(n+2)}\mid\!K_{i\bar{j}}\!\mid^{2}\\
			& \leq \mid\!B_{i\bar{j}k\bar{l}}\!\mid^{2}
			+\frac{1}{2}\mid\!\nabla_{\bar{l}}T_{ki\bar{j}}\!\mid^2,
		\end{split}
	\end{equation}
	where we have used Lemma \ref{11195} in the second equality. The above equality holds if and only if
	\[
	0=\tilde{R}_{i\bar{j}k\bar{l}}
	=\frac{1}{4}\left(\nabla_{\bar l}T_{ik\bar j}-\nabla_{\bar j}T_{ik\bar l}\right),
	\]
	where we have used \eqref{11022} in the second equality. Indeed, by definition of $\tilde R_{i\bar j k\bar l}$ and \eqref{11226},
	\begin{equation}\label{11212}
		g^{k\bar l}\tilde R_{i\bar j k\bar l}=\frac{1}{4}K_{i\bar j}.
	\end{equation}
	Hence, $\tilde{R}_{i\bar{j}k\bar{l}}=0$ implies $K_{i\bar j}=0$.
	
	For $n=3$, there is a pointwise isomorphism (e.g., \cite[Lemma 6.20]{Voisin})
	\[
	\Lambda_\omega:\Lambda^{2,2}\left(T_x^{1,0}X\right)^{\ast}\to\Lambda^{1,1}\left(T_x^{1,0}X\right)^{\ast}.
	\]
	Together with \eqref{11212}, we get
	\begin{equation}\label{11258}
		\tilde R_{i\bar j k\bar l}=\frac{1}{4}\left(K_{i\bar j}g_{k\bar l}+K_{k\bar l}g_{i\bar j}-K_{i\bar l}g_{k\bar j}-K_{k\bar j}g_{i\bar l}\right).
	\end{equation}
	Direct calculation yields
	\[
	\mid\!\tilde R_{i\bar j k\bar l}\!\mid^{2}=\frac{1}{4}\mid\!K_{i\bar j}\!\mid^{2}.
	\]
	By \eqref{11213}, we obtain
	\[
	\begin{split}
		\frac{1}{4}\,\mathcal{D}_{\nabla}
		& = \mid\!B_{i\bar{j}k\bar{l}}\!\mid^{2}
		+\frac{1}{2}\mid\!\nabla_{\bar{l}}T_{ki\bar{j}}\!\mid^2
		-\frac{6}{5}\mid\!K_{i\bar{j}}\!\mid^{2}\\
		& \leq \mid\!B_{i\bar{j}k\bar{l}}\!\mid^{2}
		+\frac{1}{2}\mid\!\nabla_{\bar{l}}T_{ki\bar{j}}\!\mid^2,
	\end{split}
	\]
	and the equality holds if and only if $K_{i\bar{j}}=0$.
	
	Suppose $(X, \chi)$ is a compact astheno-K\"{a}hler manifold. Direct calculation yields that $\tilde{R}_{i\bar{j}k\bar{l}}=0$ if and only if
	\[
	\sqrt{-1}\partial\bar{\partial}\,\omega=-\frac{1}{4}g^{p\bar{q}}T_{ki\bar{q}}\bar{T}_{lj\bar{p}}\mathrm{d}z^{k}\wedge\mathrm{d}\bar{z}^{l}\wedge\mathrm{d}z^{i}\wedge\mathrm{d}\bar{z}^{j}.
	\]
	Then we have
	\[
	\begin{split}
		0
		& = \int_X\omega\wedge\frac{\sqrt{-1}\partial\bar\partial\,\chi^{n-2}}{(n-2)!}
		=\int_{X}\sqrt{-1}\partial\bar{\partial}\,\omega\wedge\frac{\chi^{n-2}}{(n-2)!}\\
		& = \frac{1}{2}\int_{X}\chi^{i\bar{j}}\chi^{k\bar{l}}g^{p\bar{q}}T_{ki\bar{q}}\bar{T}_{lj\bar{p}}\, \frac{\chi^{n}}{n!}\geq0,
	\end{split}
	\]
	which implies $T=0$, i.e., $\omega$ is K\"{a}hler Ricci-flat.
	
	Conversely, a Ricci-flat K\"ahler metric has $T=0$, $\kappa_{\nabla}=0$, so equality holds.
\end{proof}

\begin{example}\label{11248}
	There exists a compact complex threefold $X$ admitting no pluriclosed
	metric and carrying a family of balanced metrics
	$\{\omega_\alpha\}_{\alpha>0}$. For each $\alpha>0$, the Chern
	connection $\nabla$ of $\omega_\alpha$ and its torsion $T$ satisfy
	\[
	T\ne0,\quad \nabla T=0,\quad
	\rho_\nabla=\kappa_\nabla=0,\quad
	R_\nabla=B_\nabla\not\equiv0.
	\]
	Moreover,
	\begin{equation}\label{11242}
		\mathcal N=0, \quad\mathcal{D}_{\nabla}=4\mid\!B_\nabla\!\mid^2.
	\end{equation}
	The second Bott Chern class and the weighted integral satisfy
	\[
	c_2^{\mathrm{BC}}(X)=0,\quad
	\int_Xc_2\left(X,\omega_\alpha\right)\wedge\omega_\alpha=48\alpha>0.
	\]
\end{example}
\begin{proof}
	Use the holomorphic torus bundle of Goldstein-Prokushkin
	\cite[Theorem 1]{GP} with twisting form of pure type (2,0), as allowed
	in \cite[(3.4), (5.7)]{BBFTY}. The Chern form and fibre integral
	calculations for this ansatz are in \cite[(4.2)-(4.6)]{BBFTY}.
	For the anti-self-dual (1,1) twist, see Fu-Yau
	\cite[Propositions 6, 8 and Lemma 13]{FY}.
	
	Biswas and Dumitrescu construct a principal elliptic bundle over a Kummer K3 surface using the holomorphic two-form induced by $\mathrm dz^1\wedge\mathrm dz^2$
	\cite[proof of Theorem~6.1]{BDKummer}. On an elliptic bundle of this type, Ni obtains nonzero parallel Chern torsion and computes the Chern curvature as the pullback of the curvature of a Ricci-flat K3 metric on the horizontal summand \cite[Section 10, (10.1)-(10.3)]{NiParallel}. Here we fix the fibre periods to compute the weighted second Chern integral and to exhibit equality in Theorem \ref{11191} with both torsion and curvature nonzero.
	
	For $E=\mathbb C/\left(\mathbb Z+\sqrt{-1}\mathbb Z\right)$, let $A=E\times E$, and let $S$ be the Kummer $K3$ surface obtained by resolving the sixteen singular points of $A/\{\pm1\}$.
	
	Denote the blowup at the sixteen fixed points and the double cover by
	\[
	p:\widetilde A\to A, \quad q:\widetilde A\to S,
	\]
	respectively. Normalize the holomorphic 2-form $\Omega$ on $S$ by
	\[
	q^*\Omega=2\,p^*\left(\mathrm dz^1\wedge\mathrm dz^2\right).
	\]
	
	For every integral two-cycle $C$ represented by an oriented surface transverse to the branch divisor,
	\[
	\int_C\Omega
	=\int_{q^{-1}(C)}p^*\left(\mathrm dz^1\wedge\mathrm dz^2\right)
	\in\mathbb Z+\sqrt{-1}\mathbb Z.
	\]
	Moreover,
	\[
	\int_S\Omega\wedge\bar\Omega
	=2\int_A\mathrm dz^1\wedge\mathrm dz^2
	\wedge\mathrm d\bar z^1\wedge\mathrm d\bar z^2=8.
	\]
	Choose a Ricci-flat K\"{a}hler metric $\omega_S$ with $\nabla^S$ being its Chern connection. It is well-Known that $\nabla^S\Omega=0$.
	
	Since $\operatorname{Re}\Omega$ and $\operatorname{Im}\Omega$ have integral periods, take the corresponding principal two-torus bundle $\pi:X\to S$, with real connection forms $\eta_1$, $\eta_2$ of period one. Thus
	\[
	\mathrm d\eta_1=\pi^*\operatorname{Re}\Omega,\quad
	\mathrm d\eta_2=\pi^*\operatorname{Im}\Omega.
	\]
	Set $\zeta=\eta_1+\sqrt{-1}\eta_2$. The pullbacks of the (1,0)-forms on $S$, together with $\zeta$, define an integrable complex structure. Thus
	\[
	\bar\partial\zeta=0,\quad \partial\zeta=\pi^*\Omega.
	\]
	For every $\alpha>0$, define
	\[
	\omega_\alpha=\pi^*\omega_S+\alpha\sqrt{-1}\,\zeta\wedge\bar\zeta.
	\]
	Then
	\[
	\mathrm d\zeta=\pi^*\Omega\in\mathcal A^{2,0}(X),\quad
	\mathrm d\omega_\alpha^2
	=2\alpha\sqrt{-1}\,\pi^*\omega_S\wedge
	\left(\pi^*\Omega\wedge\bar\zeta-\zeta\wedge\pi^*\bar\Omega\right)=0.
	\]
	
	Let $e^1$, $e^2$ be a local unitary coframe on $S$, with Chern connection matrix $\gamma$. In the unitary coframe $\varphi=\left(\pi^*e^1,\pi^*e^2,\sqrt{\alpha}\,\zeta\right)^t$ on \(X\), the first structure equation reads
	\[
	\Gamma=\begin{pmatrix}\pi^*\gamma&0\\0&0\end{pmatrix},\quad
	\mathrm d\varphi+\Gamma\wedge\varphi
	=\begin{pmatrix}0\\0\\\sqrt{\alpha}\,\pi^*\Omega\end{pmatrix}.
	\]
	Since $\bar{\Gamma}^{t}=-\Gamma$ and the torsion forms have type (2,0), $\Gamma$ is the Chern connection matrix of $\omega_\alpha$. Writing $V_3$ for the third vector of the dual frame, we obtain
	\[
	\nabla V_3=0,\quad
	T=\sqrt{\alpha}\,\pi^*\Omega\otimes V_3,\quad
	\nabla T=0.
	\]
	The second structure equation gives
	\[
	R_{\nabla^S}=\mathrm d\gamma+\gamma\wedge\gamma,\quad
	R_{\nabla}
	=\mathrm d\Gamma+\Gamma\wedge\Gamma
	=\begin{pmatrix}\pi^*R_{\nabla^S}&0\\0&0\end{pmatrix}.
	\]
	Hence,
	\[
	R_{i\bar j k\bar l}=
	\begin{cases}
		R^{\nabla^S}_{i\bar j k\bar l}\circ\pi,& i,j,k,l\in\{1,2\},\\
		0,&\text{if any index equals }3.
	\end{cases}
	\]
	The K\"{a}hler symmetries and Ricci-flatness of $\omega_S$ imply
	\[
	R_{i\bar j}=K_{i\bar j}=R^3_{i\bar j}=R^4_{i\bar j}=0,
	\quad R_\nabla=B_\nabla.
	\]
	Consequently, \eqref{11242} hold.
	
	On each fibre $F_y=\pi^{-1}(y)$,
	\[
	\omega_\alpha\mid_{F_y}=2\alpha\,\left(\eta_1\wedge\eta_2\right)\mid_{F_y},
	\quad \int_{F_y}\eta_1\wedge\eta_2=1.
	\]
	By Chern?Weil theory,
	\[
	c_2\left(X,\omega_\alpha\right)=\pi^*c_2\left(S,\omega_S\right).
	\]
	Therefore
	\[
	\int_Xc_2\left(X,\omega_\alpha\right)\wedge\omega_\alpha
	=2\alpha\int_Sc_2(S,\omega_S)=48\alpha>0.
	\]
	In particular, $B_\nabla\not\equiv0$.
	
	Since
	\[
	\int_S\left(c_2\left(S,\omega_S\right)-3\,\Omega\wedge\bar\Omega\right)=0, 
	\]
	there is a $f\in\mathcal{A}^{0}_{\mathbb R}(S)$ satisfying
	\[
	\sqrt{-1}\partial\bar\partial f\wedge\omega_S
	=c_2\left(S,\omega_S\right)-3\,\Omega\wedge\bar\Omega.
	\]
	By
	\[
	\bar\partial\left(\zeta\wedge\bar\zeta\right)=-\zeta\wedge\pi^*\bar\Omega,
	\]
	we obtain
	\[
	\begin{aligned}
		\sqrt{-1}\partial\bar\partial\,\omega_\alpha
		&=\alpha\,\pi^*\left(\Omega\wedge\bar\Omega\right),\\
		c_2\left(X,\omega_\alpha\right)
		&=\sqrt{-1}\partial\bar\partial
		\left(\pi^*(f\omega_S)+3\sqrt{-1}\,\zeta\wedge\bar\zeta\right).
	\end{aligned}
	\]
	Thus $c_2^{\mathrm{BC}}(X)=0$. Moreover, direct calculation yields
	\[
	\int_X\mathcal D_\nabla=1536\,\pi^2\alpha, \quad
	\int_X\mathcal L_\nabla=-1920\,\pi^2\alpha.
	\]
	
	By a similar argument after Corollary 1.3 in \cite{FLY}, $X$ admits no pluriclosed metric.
\end{proof}

\begin{proof}[Proof of Theorem \ref{11194}]
	Let $(X, \omega)$ be a connected Gauduchon surface.
	
	Fix $x\in X$. Let $S_x\subset T_x^{1,0}X$ be the unit sphere. Denote $\mathrm{d}S_x$ to be the volume induced by $\omega$ at $x$. Denote
	\[
	\mathrm{d}\sigma_x=\frac{\mathrm{d}S_x}{\int_{S_x}\mathrm{d}S_x}.
	\]
	By \cite{XKY}, the mean holomorphic sectional curvature at $x$ satisfies
	\[
	\bar{\mathcal{H}}
	=\int_{S_x}R_{i\bar j k\bar l}(x)\,v^i\bar v^jv^k\bar v^l\,\mathrm{d}\sigma_x(v)
	=\frac{1}{6}\left(s+\hat s\right).
	\]
	By \eqref{10271}, we get
	\[
	\mathcal{H}_{\min}\leq\bar{\mathcal{H}}\leq\mathcal{H}_{\max},\quad
	s=3\,\bar{\mathcal{H}}+\frac{1}{2}\mid\!\tau\!\mid^{2}\in\mathcal{I}_{x}.
	\]
	Consequently,
	\[
	\operatorname{dist}\left(0, \mathcal{I}_x\right)^2
	\leq s^2\leq\mathcal{N},
	\]
	which together with Theorem \ref{11181} implies \eqref{11232}.
	
	Suppose equality in (1.7) holds at every point of $X$. Then we have
	\[
	s^2=\mathcal{N},\quad \tilde R_{i\bar j k\bar l}=M_{i\bar j k\bar l}=0,\quad R^\circ_{i\bar j}=K^\circ_{i\bar j}=0.
	\]
	By Proposition \ref{11221}, we get
	\[
	R_{i\bar j k\bar l}=\hat R_{i\bar j k\bar l},
	\]
	which implies $s=\hat s$. Thus \eqref{10271} gives $\tau=0$, which together with \eqref{06041} implies $\omega$ is K\"ahler.
	
	Since $R^\circ_{i\bar j}=0$, the same arguments for the equality in the proof of Theorem \ref{11224}, we get $\omega$ is K\"ahler-Einstein, i.e., $s$ is a constant.
	
	If $s=0$, then $\omega$ is K\"ahler Ricci-flat. If $s>0$, we get
	\[
	\operatorname{dist}\left(0, \mathcal{I}_x\right)=s\in\mathcal{I}_x,
	\]
	which implies $s=3\,\mathcal{H}_{\min}$. Since $s=3\,\bar{\mathcal H}$ a prior, the continuous nonnegative function
	\[
	h_x-\mathcal H_{\min}(x)
	\]
	has mean zero on $S_x$. If $s<0$, the same identity gives $s=3\,\mathcal H_{\max}$, and one applies the same argument to $\mathcal H_{\max}(x)-h_x$. In both cases, for every $v\in S_x$,
	\[
	h_x(v)=\frac{s}{3}.
	\]
	Thus $\omega$ has nonzero constant holomorphic sectional curvature.
	
	Conversely, for a Ricci-flat K\"ahler metric, we have
	\[
	\tau=0, \quad \mathcal{N}=0, \quad \bar{\mathcal{H}}=0\in\mathcal{I}_x.
	\]
	For a K\"ahler metric of constant holomorphic sectional curvature $c\neq0$, we have
	\[
	R_{i\bar j k\bar l}
	=\frac c2\left(g_{i\bar j}g_{k\bar l}+g_{i\bar l}g_{k\bar j}\right),
	\]
	which implies
	\[
	s=3\,c,\quad \mathcal{I}_x=\{3\,c\},\quad \mathcal{N}=s^2
	\]
	Since $\tau=0$, equality in \eqref{11232} holds in both cases.
\end{proof}

\section{Proof of Theorem \ref{11243}}\label{112410}

For balanced metrics, we estimate $\mathcal D_\nabla$ and $\mathcal L_\nabla$ and derive a Ricci pinching condition that ensures $\mathcal D_\nabla\geq0$. On a compact balanced threefold admitting a pluriclosed metric, we show that the given metric is K\"ahler precisely when its first and second Chern-Ricci forms coincide. We also give a pointwise formula for $\mid\!B_\nabla\!\mid^2$ in terms of the holomorphic sectional curvature of an arbitrary Hermitian metric.

\begin{proof}[Proof of Theorem \ref{11243}]
	Let $(X, \omega)$ be a balanced manifold. By \eqref{11022} and \eqref{10271}, we have $s=\hat{s}$ and
	\[
	\hat{R}_{i\bar j}=\frac14\left(3R_{i\bar j}+K_{i\bar j}\right),
	\]
	which implies
	\begin{equation}\label{112411}
		\frac34\mid\!K_{i\bar j}-R_{i\bar j}\!\mid^2
		=3\mid\!R_{i\bar j}^\circ\!\mid^2+\mid\!K_{i\bar j}^\circ\!\mid^2
		-4\mid\!\hat R_{i\bar j}^\circ\!\mid^2.
	\end{equation}
	
	By Corollary 1 in \cite{Yang},
	\[
	\mid\!\hat{R}_{i\bar jk\bar l}\!\mid^{2}
	=\mid\!B_{i\bar jk\bar l}\!\mid^{2}
	+\frac{4}{n+2}\mid\!\hat{R}_{i\bar{j}}\!\mid^{2}-\frac{2s^{2}}{(n+1)(n+2)}.
	\]
	Combining this with \eqref{11227} yields
	\begin{equation}\label{112412}
		\mid\!R_{i\bar jk\bar l}\!\mid^{2}
		=\mid\!B_{i\bar jk\bar l}\!\mid^{2}
		+\mid\!\tilde{R}_{i\bar jk\bar l}\!\mid^{2}+\mid\!M_{i\bar jk\bar l}\!\mid^{2}
		+\frac4{n+2}\mid\!\hat R_{i\bar j}^\circ\!\mid^2+\frac{2s^2}{n(n+1)}.
	\end{equation}
	Inserting \eqref{112411} and \eqref{112412} into \eqref{11222} yileds
	\begin{equation}\label{112413}
		\begin{split}
			\mathcal L_\nabla
			& = -(n+2)\left(\mid\!B_{i\bar jk\bar l}\!\mid^{2}+\mid\!\tilde{R}_{i\bar jk\bar l}\!\mid^{2}+\mid\!M_{i\bar jk\bar l}\!\mid^{2}\right)\\
			& \quad
			+\frac34\mid\!K_{i\bar j}-R_{i\bar j}\!\mid^2
			+(n+1)\left(\mid\!K_{i\bar j}^\circ\!\mid^2-\mid\!R^\circ_{i\bar j}\!\mid^2\right)+\frac{n(n-1)}{n+1}s^2.
		\end{split}
	\end{equation}
	For the case $\rho_\nabla=\kappa_\nabla$, combining Proposition \ref{11221} and \eqref{112413} yields \eqref{11253}.
	
	Since
	\[
	g^{i\bar j}g^{k\bar l}\tilde R_{i\bar j k\bar l}
	=\frac{1}{4}\,g^{i\bar j}\left(K_{i\bar j}-R_{i\bar j}\right)=0,
	\]
	we calculate
	\begin{equation}\label{11251}
		\begin{split}
			\mid\!\tilde{R}_{i\bar jk\bar l}\!\mid^{2}
			-\frac{1}{4(n-2)}\mid\!K_{i\bar j}-R_{i\bar j}\!\mid^2
			& = \frac{1}{16(n-2)^2}\,\Big|4(n-2)\tilde R_{i\bar j k\bar l}\\
			& \quad
			-\left(K_{i\bar j}-R_{i\bar j}\right)g_{k\bar l}
			-\left(K_{k\bar l}-R_{k\bar l}\right)g_{i\bar j}\\
			& \quad
			+\left(K_{i\bar l}-R_{i\bar l}\right)g_{k\bar j}
			+\left(K_{k\bar j}-R_{k\bar j}\right)g_{i\bar l}\Big|^{2}\\
			& \geq0.
		\end{split}
	\end{equation}
	Since
	\[
	-g^{k\bar l}M_{k\bar l i\bar j}
	=g^{k\bar l}M_{i\bar j k\bar l}
	=\frac12(K_{i\bar j}-R_{i\bar j}),
	\]
	we calculate
	\begin{equation}\label{11252}
		\begin{split}
			\mid\!M_{i\bar jk\bar l}\!\mid^{2}
			-\frac{1}{2n}\mid\!K_{i\bar j}-R_{i\bar j}\!\mid^2
			& = \frac{1}{4n^2}\,\Big|2n\,M_{i\bar j k\bar l}
			-\left(K_{i\bar j}-R_{i\bar j}\right)g_{k\bar l}
			+\left(K_{k\bar l}-R_{k\bar l}\right)g_{i\bar j}\Big|^{2}\\
			& \geq0.
		\end{split}
	\end{equation}
	By \eqref{11251} and \eqref{11252}, we obtain
	\[
	(n-2)\mid\!\tilde R_{i\bar j k\bar l}\!\mid^2+n\mid\!M_{i\bar j k\bar l}\!\mid^2
	\geq\frac34\mid\!K_{i\bar j}-R_{i\bar j}\!\mid^2.
	\]
	Hence
	\[
	\begin{split}
		(n+2)\left(\mid\!\tilde R_{i\bar j k\bar l}\!\mid^2
		+\mid\!M_{i\bar j k\bar l}\!\mid^2\right)
		-\frac34\mid\!K_{i\bar j}-R_{i\bar j}\!\mid^2
		& \geq 4\mid\!\tilde R_{i\bar j k\bar l}\!\mid^2 +2\mid\!M_{i\bar j k\bar l}\!\mid^2=\mid\!\nabla''T\!\mid^2 \\
		& \geq \frac{2(n-1)}{n(n-2)}\mid\!K_{i\bar j}-R_{i\bar j}\!\mid^2,
	\end{split}
	\]
	where the equality follows from Lemma \ref{11195}, and the last inequality follows from \eqref{11251} and \eqref{11252}. Using the bound by
	$\tfrac{2(n-1)}{n(n-2)}\mid\!K_{i\bar j}-R_{i\bar j}\!\mid^2$ in \eqref{112413} and applying Proposition \ref{11246} gives \eqref{11244}. The bound by $\mid\!\nabla''T\!\mid^2$ in \eqref{112413} gives \eqref{11245}. 
\end{proof}

\begin{corollary}\label{11254}
	Let $(X,\omega)$ be a connected balanced manifold. If
	\begin{equation}\label{11257}
		(n+1)^2\mid\!\kappa^\circ\!\mid^2\leq(2n+1)\mid\!\rho^\circ\!\mid^2,
	\end{equation}
	then $\mathcal D_\nabla\ge0$. At each point satisfying \eqref{11257}, $\mathcal D_\nabla=0$ if and only if
	\[
	R_{i\bar j k\bar l}=\frac{s}{n(n+1)}
	\left(g_{i\bar j}g_{k\bar l}+g_{i\bar l}g_{k\bar j}\right).
	\]
	If \eqref{11257} holds on $X$, then $\mathcal D_\nabla\equiv0$ if and only if $\omega$ is Chern-flat or is K\"ahler with nonzero constant holomorphic sectional curvature.
\end{corollary}
\begin{proof}
	For $n\geq3$, by \eqref{11251} and \eqref{11252}, we get
	\begin{equation}\label{11256}
		\frac34\mid\!K_{i\bar j}-R_{i\bar j}\!\mid^2
		\le(n-2)\mid\!\tilde R_{i\bar j k\bar l}\!\mid^2+n\mid\!M_{i\bar j k\bar l}\!\mid^2,
	\end{equation}
	which is also true for $n=2$ since balancedness is equivalent to Kahlerness.
	
	By Proposition \ref{11246}, Lemma \ref{11195}, \eqref{112413} and \eqref{11256},
	\[
	\begin{split}
		(n+2)\mathcal D_\nabla
		&\ge (n+1)(n+2)\mid\!B_\nabla\!\mid^2+(2n+1)\mid\!\rho_\nabla^\circ\!\mid^2
		-(n+1)^2\mid\!\kappa_\nabla^\circ\!\mid^2\\
		&\quad +(n+1)\mid\!\nabla''T\!\mid^2.
	\end{split}
	\]
	At a point satisfying \eqref{11257}, the right hand side of the above identity is nonnegative. If $\mathcal D_\nabla=0$ there, then
	\[
	B_\nabla=0,\quad \nabla''T=0.
	\]
	Lemma \ref{11195} gives
	\[
	\tilde R_{i\bar j k\bar l}=M_{i\bar j k\bar l}=0,
	\]
	which implies $K_{i\bar j}=R_{i\bar j}$. Now \eqref{11257} gives $\rho_\nabla^\circ=0$. Thus $R_\nabla=\hat R_\nabla$, and Definition \ref{09131} gives the stated curvature formula. Its converse follows from \eqref{11222}.
	
	If this formula holds throughout $X$, then $\rho_\nabla=\frac{1}{n}s\,\omega$
	and $\mathrm d\left(s\omega\right)=0$. Balancedness gives $\mathrm ds\wedge\omega^{n-1}=0$. Wedge multiplication by $\omega^{n-1}$ is injective on one-forms, so $s$ is constant on $X$. For $s=0$, the curvature vanishes. For $s\ne0$, we have $\mathrm d\omega=0$ and
	\[
	h_x(v)=\frac{2s}{n(n+1)}
	\]
	for every $v\in S_x$. Conversely, both classes of metrics satisfy the stated curvature formula.
\end{proof}

\begin{proposition}\label{11249}
	Let $(X,\omega)$ be a balanced threefold. Thefollowing identity holds
	\[
	4\mid\!\tilde R_{i\bar j k\bar l}\!\mid^2=\mid\!K_{i\bar j}-R_{i\bar j}\!\mid^2.
	\]
	At each point where $\rho_\nabla=\kappa_\nabla$,
	\[
	\frac32\,s^2-\mathcal L_\nabla
	=5\mid\!B_\nabla\!\mid^2+\frac52\mid\!\nabla''T\!\mid^2.
	\]
	If $X$ is compact and admits a pluriclosed metric $\chi$, then $\omega$ is K\"{a}hler if and only if
	\[
	\int_X\left(\kappa_\nabla-\rho_\nabla\right)\wedge\omega\wedge\chi=0,
	\]
	which is equivalent to $\rho_\nabla=\kappa_\nabla$ on $X$.
\end{proposition}
\begin{proof}
	Balancedness and \eqref{11022} give
	\[
	R^3_{i\bar j}=R^4_{i\bar j}=R_{i\bar j}.
	\]
	Hence
	\[
	g^{i\bar j}g^{k\bar l}\tilde R_{i\bar j k\bar l}
	=\frac{1}{4}g^{i\bar j}\left(K_{i\bar j}-R_{i\bar j}\right)=0.
	\]
	On threefold, the pointwise Lefschetz contraction $\Lambda_\omega:\Lambda^{2,2}\to\Lambda^{1,1}$ is injective. The tensors on the two sides of the following identity alternate in $i,k$ and in $j,l$, and their contractions with $g^{k\bar l}$ agree. It follows that
	\begin{equation}\label{11262}
		\tilde R_{i\bar j k\bar l}
		=\frac{1}{4}\left(
		\left(K_{i\bar j}-R_{i\bar j}\right)g_{k\bar l}
		+\left(K_{k\bar l}-R_{k\bar l}\right)g_{i\bar j}
		-\left(K_{i\bar l}-R_{i\bar l}\right)g_{k\bar j}
		-\left(K_{k\bar j}-R_{k\bar j}\right)g_{i\bar l}\right),
	\end{equation}
	which implies
	\[
	\mid\!\tilde R_{i\bar j k\bar l}\!\mid^2
	=\frac{1}{4}\mid\!K_{i\bar j}-R_{i\bar j}\!\mid^2.
	\]
	When $\rho_\nabla=\kappa_\nabla$, we have $\tilde R=0$. By Lemma \ref{11195},
	\[
	2\mid\!M_{i\bar j k\bar l}\!\mid^2=\mid\!\nabla_{\bar{l}}T_{ki\bar j}\!\mid^2,
	\]
	and \eqref{11253} gives the desired formula for $\mathcal L_\nabla$.
	
	Direct calculation yields
	\[
	\sqrt{-1}\,\partial\bar\partial\,\omega
	=
	\left(\tilde R_{i\bar j k\bar l}-\frac14\,g^{p\bar q}T_{ki\bar q}\bar{T}_{lj\bar p}\right)\mathrm dz^k\wedge\mathrm d\bar z^l
	\wedge\mathrm dz^i\wedge\mathrm d\bar z^j.
	\]
	The four summands of $4\widetilde R$ give the same wedge product
	in the first line. Since $\kappa_\nabla-\rho_\nabla
	=\sqrt{-1}(K_{i\bar j}-R_{i\bar j})\mathrm dz^i\wedge\mathrm d\bar z^j$,
	that line equals $(\rho_\nabla-\kappa_\nabla)\wedge\omega$.
	
	Assume that $X$ is compact and $\partial\bar\partial\chi=0$. Using \eqref{11262}, we calculate
	\[
	\begin{aligned}
		0&=\int_X\sqrt{-1}\,\omega\wedge\partial\bar\partial\,\chi
		=\int_X\sqrt{-1}\,\partial\bar\partial\,\omega\wedge\chi\\
		&=\int_X\left(\rho_\nabla-\kappa_\nabla\right)\wedge\omega\wedge\chi
		+\frac12\int_X\chi^{i\bar j}\chi^{k\bar l}g^{p\bar q}T_{ki\bar q}\bar{T}_{lj\bar p},
	\end{aligned}
	\]
	where the second integrand is nonnegative and vanishes exactly when $T=0$. Consequently the stated integral vanishes precisely when $\mathrm d\omega=0$. For a K\"ahler metric the curvature symmetries give $\rho_\nabla=\kappa_\nabla$, which completes the equivalences.
\end{proof}
Fix a pluriclosed metric $\chi$ on a compact complex threefold $X$.
By Proposition~\ref{11249}, $X$ admits a K\"ahler metric if and only
if there exists a balanced metric $\omega$ such that
\[
\int_X(\kappa_\nabla-\rho_\nabla)\wedge\omega\wedge\chi=0.
\]
In the forward direction one may take $\omega$ to be K\"ahler. Thus, in dimension three, the Fino-Vezzoni conjecture asks whether such an $\omega$ can be chosen whenever $X$ admits a balanced metric.
\begin{remark}
	On a balanced threefold, Proposition \ref{11249} gives
	\[
	4\mid\!\tilde R\!\mid^2=\mid\!\kappa_\nabla-\rho_\nabla\!\mid^2.
	\]
	Hence the sharp lower bound in Theorem \ref{11181} is attained at $x$ if and only if $\rho_\nabla(x)=\kappa_\nabla(x)=0$.
	
	In any dimension, suppose that $X$ is compact and
	\[
	\mathcal N=\frac{1}{2}\,(n+1)\mid\!\nabla''T\!\mid^2.
	\]
	Then we get
	\[
	\tilde R=0, \quad \rho_\nabla=\kappa_\nabla=0.
	\]
	Thus we obtain
	\[
	g^{i\bar j}g^{k\bar l}\tilde R_{i\bar j k\bar l}
	=\frac{1}{2}\left(s-\hat s\right)=0,
	\]
	which together with \eqref{11261} implies $\omega$ is balanced. If $X$ also admits an astheno-K\"ahler metric, then $\tilde R=0$ gives \eqref{11235}, and Theorem \ref{11234} implies that $\omega$ is K\"ahler Ricci-flat. Conversely, every Ricci-flat K\"ahler metric attains the sharp bound.
\end{remark}

Du related Chern numbers of compact K\"ahler-Einstein manifolds to integrals of the holomorphic sectional curvature \cite{Du}. For an arbitrary Hermitian metric, the squared norm of the Bochner curvature in Definition \ref{09131} has the following pointwise integral representation.

\begin{proposition}\label{11259}
	Let $(X,\omega)$ be a Hermitian manifold of complex dimension $n$.
	With the notation of Theorem \ref{11194}, for every $x\in X$, $v\in S_x$, set
	\[
	h_{0,x}(v)=h_x(v)-\frac{1}{n(n+1)}\left(s+\hat s\right)(x)
	-\frac{4}{n+2}\hat R^\circ_{i\bar j}(x)v^i\bar v^j
	\]
	The following identity holds
	\[
	\mid\!B_\nabla\!\mid^2(x)
	=\frac{1}{4}\,n(n+1)(n+2)(n+3)\,
	\int_{S_x}\mid\!h_{0,x}(v)\!\mid^2\,\mathrm{d}\sigma_x(v).
	\]
\end{proposition}
\begin{proof}
	Fix $x\in X$. For $v\in S_x$, we obtain from Definition \ref{09131} that
	\[
	h_{0,x}(v)=B_{i\bar j k\bar l}(x)v^i\bar v^jv^k\bar v^l.
	\]
	The tensor $B_\nabla(x)$ has the K\"ahler curvature symmetries
	and all its traces vanish. Thus \cite[Corollary~4]{Magnusson},
	applied to $B_\nabla(x)$, gives
	\[
	\int_{S_x}\mid\!h_{0,x}(v)\!\mid^2\,\mathrm{d}\sigma_x(v)
	=\frac{1}{\binom{n+1}{2}\binom{n+3}{2}}\mid\!B_\nabla\!\mid^2(x)
	\]
	as claimed.
\end{proof}

\end{document}